\documentclass{elsart}
\usepackage{amssymb}
\usepackage{amsfonts}
\usepackage{amsmath}

\input{tcilatex}
\begin{document}

\begin{frontmatter}%

\title{The explicit formula of Hankel determinant with Catalan elements}%

\author[1]{Jishe Feng}
\thanks[1]{Corresponding author. E-mail: gsfjs6567@126.com.}%

\address{Department of Mathematics, Longdong University,  Qingyang,  Gansu,  745000,  China
E-mail: gsfjs6567@126.com.}%

\begin{abstract}%

Applying Johann Cigler's Hankel determinant formula in terms of the binomial
coefficient determinants, which is simplified from Christian Krattenthale's,
we get an explicit formula of Hankel determinants for general. As far as I
know, those are new results.

2010 Mathematics Subject Classification: 05A15, 11B83.

\end{abstract}%

\begin{keyword}%
Hankel determinant; Hankel transform; Catalan number; binomial coefficient%
\end{keyword}%

\end{frontmatter}%

\section{Introduction}

The Catalan numbers $C_{n}$ are the famous combinatorial numbers. Stanley 
\cite{Stanley} lists 214 different kinds of objects counted by Catalan
numbers. Many researchers discuss the determinant formulas of matrices with
Catalan elements or its representation in terms of determinants \cite%
{Kreweras}. There are a lot of relations between the determinant of the
Hessenberg matrix and many well-known number sequences (see \cite{Feng} \cite%
{Yilmaz} and references therein).

The research in Hankel matrices is one of the primary topics in
combinatorial matrices. The Hankel determinants evaluation has attracted
attention since the 19th century. The \textit{Hankel transform} of the
sequence\ $\{a_{0},a_{1},a_{2},...\}$ is the Hankel determinant of the form $%
det(a_{i+j+r})_{i,j=0}^{n-1}$, where $r$ is a non-negative integer (see \cite%
{Layman}). For $r=0,1,2$, and let $(a_{n})_{n\geq 0}$ be the Catalan squence 
$(C_{n})_{n\geq 0}$, there are many results \cite{Cigler}\cite{Chammam}\cite%
{Kreweras}\cite{kk}\cite{Wang}. But there is no result for $n>5$.

In this paper, we will show how to get the following Hankel determinants
formulas:%
\begin{equation}
\det_{0\leq i,j\leq n-1}(C_{i+j+4})=\frac{4}{3!\cdot 5!}%
(n+1)(n+2)^{2}(n+3)(2n+3)(2n+5),  \label{cc1}
\end{equation}%
\begin{eqnarray}
\det_{0\leq i,j\leq n-1}(C_{i+j+5}) &=&\frac{8}{5!\cdot 7!}%
(n+1)(n+2)^{2}(n+3)^{2}(n+4)  \label{cc2} \\
&&(2n+3)(2n+5)^{2}(2n+7),  \notag
\end{eqnarray}%
\begin{eqnarray}
\det_{0\leq i,j\leq n-1}(C_{i+j+6}) &=&\frac{2^{5}}{5\cdot 7!\cdot 9!}%
(n+1)(n+2)^{2}(n+3)^{3}(n+4)^{2}(n+5)  \label{cc3} \\
&&(2n+3)(2n+5)^{2}(2n+7)^{2}(2n+9),  \notag
\end{eqnarray}%
\begin{eqnarray}
\det_{0\leq i,j\leq n-1}(C_{i+j+7}) &=&\frac{3\cdot 2^{10}}{7!\cdot 9!\cdot
11!}(n+1)(n+2)^{2}(n+3)^{3}(n+4)^{3}(n+5)^{2}  \label{cc4} \\
&&(n+6)(2n+3)(2n+5)^{2}(2n+7)^{3}(2n+9)^{2}(2n+11),  \notag
\end{eqnarray}%
and general caese%
\begin{eqnarray}
\det_{0\leq i,j\leq n-1}(C_{i+j+r}) &=&\frac{(n+1)(n+2)\cdots (n+r-1)}{%
3!\cdot 5!\cdot 7!\cdots (2r-3)!}  \label{d2} \\
&&\cdot (2n+3)(2n+5)\cdots (2n+2r-3)\cdot (r-2)!\newline
\notag \\
&&\cdot (2n+4)(2n+6)\cdots (2n+2r-4)\cdot (r-3)!\newline
\notag \\
&&\cdot (2n+5)(2n+7)\cdots (2n+2r-5)\cdot (r-4)!\newline
\   \notag \\
&&\cdot (2n+6)(2n+8)\cdots (2n+2r-6)\cdot (r-5)!\newline
\notag \\
&&\cdots \cdots \newline
\ \ \ \ \ \ \ \ \ \ \ \ \ \ \ \ \ \ \   \notag \\
&&\cdot (2n+r-1)(2n+r+1)\cdot 2!  \notag \\
&&\cdot (2n+r)\cdot 1!.  \notag
\end{eqnarray}

\section{The determinants of Hankel and binomial coefficients}

In \cite{Cigler} Cigler simplify Christian Krattenthaler's result \cite{kk}
about the binomial coefficients and Hankel\ determinants as the follow,%
\begin{equation}
\det_{0\leq i,j\leq n-1}(C_{i+j+r})=\det_{0\leq i,j\leq r-1}(\binom{i+j+n}{%
i-j+n}).  \label{d1}
\end{equation}

For $r=4$, applying this formula, we perform the $4\times 4$ determinant by
elementary operations and Laplace expansion, get a $3\times 3$ determinant,
recurrently to $2\times 2$ determinant. Finally, we obtain the result.%
\newline
$\det_{0\leq i,j\leq n-1}(C_{i+j+4})=$ \newline
$\left\vert 
\begin{array}{cccc}
\binom{n}{n} & \binom{n+1}{n-1} & \binom{n+2}{n-2} & \binom{n+3}{n-3} \\ 
\binom{n+1}{n+1} & \binom{n+2}{n} & \binom{n+3}{n-1} & \binom{n+4}{n-2} \\ 
\binom{n+2}{n+2} & \binom{n+3}{n+1} & \binom{n+4}{n} & \binom{n+5}{n-1} \\ 
\binom{n+3}{n+3} & \binom{n+4}{n+2} & \binom{n+5}{n+1} & \binom{n+6}{n}%
\end{array}%
\right\vert =\left\vert 
\begin{array}{cccc}
1 & \frac{(n+1)n}{2!} & \frac{(n+2)(n+1)n(n-1)}{4!} & \frac{%
(n+3)(n+2)(n+1)n(n-1)(n-2)}{6!} \\ 
1 & \frac{(n+2)(n+1)}{2!} & \frac{(n+3)(n+2)(n+1)n}{4!} & \frac{%
(n+4)(n+3)(n+2)(n+1)n(n-1)}{6!} \\ 
1 & \frac{(n+3)(n+2)}{2!} & \frac{(n+4)(n+3)(n+2)(n+1)}{4!} & \frac{%
(n+5)(n+4)(n+3)(n+2)(n+1)n}{6!} \\ 
1 & \frac{(n+4)(n+3)}{2!} & \frac{(n+5)(n+4)(n+3)(n+2)}{4!} & \frac{%
(n+6)(n+5)(n+4)(n+3)(n+2)(n+1)}{6!}%
\end{array}%
\right\vert $\newline
$=\frac{1}{2!\cdot 4!\cdot 6!}\left\vert 
\begin{array}{cccc}
1 & (n+1)n & (n+2)(n+1)n(n-1) & (n+3)(n+2)(n+1)n(n-1)(n-2) \\ 
0 & 2(n+1) & 4(n+2)(n+1)n & 6(n+3)(n+2)(n+1)n(n-1) \\ 
0 & 2(n+2) & 4(n+3)(n+2)(n+1) & 6(n+4)(n+3)(n+2)(n+1)n \\ 
0 & 2(n+3) & 4(n+4)(n+3)(n+2) & 6(n+5)(n+4)(n+3)(n+2)(n+1)%
\end{array}%
\right\vert $\newline
$=\frac{1}{2!\cdot 4!\cdot 6!}\left\vert 
\begin{array}{ccc}
2(n+1) & 4(n+2)(n+1)n & 6(n+3)(n+2)(n+1)n(n-1) \\ 
2(n+2) & 4(n+3)(n+2)(n+1) & 6(n+4)(n+3)(n+2)(n+1)n \\ 
2(n+3) & 4(n+4)(n+3)(n+2) & 6(n+5)(n+4)(n+3)(n+2)(n+1)%
\end{array}%
\right\vert $\newline
$=\frac{2\cdot 4\cdot 6}{2!\cdot 4!\cdot 6!}\left\vert 
\begin{array}{ccc}
n+1 & (n+2)(n+1)n & (n+3)(n+2)(n+1)n(n-1) \\ 
n+2 & (n+3)(n+2)(n+1) & (n+4)(n+3)(n+2)(n+1)n \\ 
n+3 & (n+4)(n+3)(n+2) & (n+5)(n+4)(n+3)(n+2)(n+1)%
\end{array}%
\right\vert $\newline
$=\frac{1}{3!\cdot 5!}(n+1)(n+2)(n+3)\left\vert 
\begin{array}{ccc}
1 & (n+2)n & (n+3)(n+2)n(n-1) \\ 
1 & (n+3)(n+1) & (n+4)(n+3)(n+1)n \\ 
1 & (n+4)(n+2) & (n+5)(n+4)(n+2)(n+1)%
\end{array}%
\right\vert $\newline
$=\frac{1}{3!\cdot 5!}(n+1)(n+2)(n+3)\left\vert 
\begin{array}{ccc}
1 & (n+2)n & (n+3)(n+2)n(n-1) \\ 
0 & 2n+3 & (4n+6)(n+3)n \\ 
0 & 2n+5 & (4n+10)(n+4)(n+1)%
\end{array}%
\right\vert $\newline
$=\frac{1}{3!\cdot 5!}(n+1)(n+2)(n+3)\left\vert 
\begin{array}{cc}
2n+3 & (4n+6)(n+3)n \\ 
2n+5 & (4n+10)(n+4)(n+1)%
\end{array}%
\right\vert $\newline
$=\frac{1}{3!\cdot 5!}(n+1)(n+2)(n+3)(2n+3)(2n+5)\left\vert 
\begin{array}{cc}
1 & 2(n+3)n \\ 
1 & 2(n+4)(n+1)%
\end{array}%
\right\vert $\newline
$=\frac{4}{3!\cdot 5!}(n+1)(n+2)^{2}(n+3)(2n+3)(2n+5).$

For $r=5$, applying this formula, we perform the $5\times 5$ determinant by
elementary operations and Laplace expansion, get a $4\times 4$ determinant,
recurrently to $2\times 2$ determinant. Finally, we obtain the result.

$\det_{0\leq i,j\leq n-1}(C_{i+j+5})=$ 

$\left\vert 
\begin{array}{ccccc}
\binom{n}{n} & \binom{n+1}{n-1} & \binom{n+2}{n-2} & \binom{n+3}{n-3} & 
\binom{n+4}{n-4} \\ 
\binom{n+1}{n+1} & \binom{n+2}{n} & \binom{n+3}{n-1} & \binom{n+4}{n-2} & 
\binom{n+5}{n-3} \\ 
\binom{n+2}{n+2} & \binom{n+3}{n+1} & \binom{n+4}{n} & \binom{n+5}{n-1} & 
\binom{n+6}{n-2} \\ 
\binom{n+3}{n+3} & \binom{n+4}{n+2} & \binom{n+5}{n+1} & \binom{n+6}{n} & 
\binom{n+7}{n-1} \\ 
\binom{n+4}{n+4} & \binom{n+5}{n+3} & \binom{n+6}{n+2} & \binom{n+7}{n+1} & 
\binom{n+8}{n}%
\end{array}%
\right\vert $\newline
$=\left\vert 
\begin{array}{ccccc}
1 & \frac{(n+1)n}{2!} & \frac{(n+2)(n+1)n(n-1)}{4!} & \frac{%
(n+3)(n+2)(n+1)n(n-1)(n-2)}{6!} & \frac{(n+4)(n+3)(n+2)(n+1)n(n-1)(n-2)(n-3)%
}{8!} \\ 
1 & \frac{(n+2)(n+1)}{2!} & \frac{(n+3)(n+2)(n+1)n}{4!} & \frac{%
(n+4)(n+3)(n+2)(n+1)n(n-1)}{6!} & \frac{(n+5)(n+4)(n+3)(n+2)(n+1)n(n-1)(n-2)%
}{8!} \\ 
1 & \frac{(n+3)(n+2)}{2!} & \frac{(n+4)(n+3)(n+2)(n+1)}{4!} & \frac{%
(n+5)(n+4)(n+3)(n+2)(n+1)n}{6!} & \frac{(n+6)(n+5)(n+4)(n+3)(n+2)(n+1)n(n-1)%
}{8!} \\ 
1 & \frac{(n+4)(n+3)}{2!} & \frac{(n+5)(n+4)(n+3)(n+2)}{4!} & \frac{%
(n+6)(n+5)(n+4)(n+3)(n+2)(n+1)}{6!} & \frac{%
(n+7)(n+6)(n+5)(n+4)(n+3)(n+2)(n+1)n}{8!} \\ 
1 & \frac{(n+5)(n+4)}{2!} & \frac{(n+6)(n+5)(n+4)(n+3)}{4!} & \frac{%
(n+7)(n+6)(n+5)(n+4)(n+3)(n+2)}{6!} & \frac{%
(n+8)(n+7)(n+6)(n+5)(n+4)(n+3)(n+2)(n+1)}{8!}%
\end{array}%
\right\vert $\newline

=$\frac{2\cdot 4\cdot 6\cdot 8}{2!\cdot 4!\cdot 6!\cdot 8!}\left\vert 
\begin{array}{ccc}
1 & (n+1)n & (n+2)(n+1)n(n-1) \\ 
0 & n+1 & (n+2)(n+1)n \\ 
0 & n+2 & (n+3)(n+2)(n+1) \\ 
0 & n+3 & (n+4)(n+3)(n+2) \\ 
0 & n+4 & (n+5)(n+4)(n+3)%
\end{array}%
\right. $

$\left. 
\begin{array}{cc}
(n+3)(n+2)(n+1)n(n-1)(n-2) & (n+4)(n+3)(n+2)(n+1)n(n-1)(n-2)(n-3) \\ 
(n+3)(n+2)(n+1)n(n-1) & (n+4)(n+3)(n+2)(n+1)n(n-1)(n-2) \\ 
(n+4)(n+3)(n+2)(n+1)n & (n+5)(n+4)(n+3)(n+2)(n+1)n(n-1) \\ 
(n+5)(n+4)(n+3)(n+2)(n+1) & (n+6)(n+5)(n+4)(n+3)(n+2)(n+1)n \\ 
(n+6)(n+5)(n+4)(n+3)(n+2) & (n+7)(n+6)(n+5)(n+4)(n+3)(n+2)(n+1)%
\end{array}%
\right\vert $\newline
$=\frac{1}{3!\cdot 5!\cdot 7!}\left\vert 
\begin{array}{ccc}
n+1 & (n+2)(n+1)n & (n+3)(n+2)(n+1)n(n-1) \\ 
n+2 & (n+3)(n+2)(n+1) & (n+4)(n+3)(n+2)(n+1)n \\ 
n+3 & (n+4)(n+3)(n+2) & (n+5)(n+4)(n+3)(n+2)(n+1) \\ 
n+4 & (n+5)(n+4)(n+3) & (n+6)(n+5)(n+4)(n+3)(n+2)%
\end{array}%
\right. $\newline
$\left. 
\begin{array}{c}
(n+4)(n+3)(n+2)(n+1)n(n-1)(n-2) \\ 
(n+5)(n+4)(n+3)(n+2)(n+1)n(n-1) \\ 
(n+6)(n+5)(n+4)(n+3)(n+2)(n+1)n \\ 
(n+7)(n+6)(n+5)(n+4)(n+3)(n+2)(n+1)%
\end{array}%
\right\vert $\newline
$=\frac{(n+1)(n+2)(n+3)(n+4)}{3!\cdot 5!\cdot 7!}\left\vert 
\begin{array}{ccc}
1 & (n+2)n & (n+3)(n+2)n(n-1) \\ 
1 & (n+3)(n+1) & (n+4)(n+3)(n+1)n \\ 
1 & (n+4)(n+2) & (n+5)(n+4)(n+2)(n+1) \\ 
1 & (n+5)(n+3) & (n+6)(n+5)(n+3)(n+2)%
\end{array}%
\right. $\newline
$\left. 
\begin{array}{c}
(n+4)(n+3)(n+2)n(n-1)(n-2) \\ 
(n+5)(n+4)(n+3)(n+1)n(n-1) \\ 
(n+6)(n+5)(n+4)(n+2)(n+1)n \\ 
(n+7)(n+6)(n+5)(n+3)(n+2)(n+1)%
\end{array}%
\right\vert $\newline
$=\cdots =\frac{8}{5!\cdot 7!}%
(n+1)(n+2)^{2}(n+3)^{2}(n+4)(2n+3)(2n+5)^{2}(2n+7).$

For the general case, we perform elementary operations and Laplace
expansions on the $r\times r$ determinant, get an ($r-1)\times (r-1)$
determinant, and so on, and so on until we get to the general explicit
formula (\ref{d2}).\newline
$\det_{0\leq i,j\leq n-1}(C_{i+j+r})=\left\vert 
\begin{array}{cccc}
\binom{n}{n} & \binom{n+1}{n-1} & \cdots  & \binom{n+r-1}{n-r+1} \\ 
\binom{n+1}{n+1} & \binom{n+2}{n} & \cdots  & \binom{n+r}{n-r+2} \\ 
\cdots  & \cdots  & \cdots  & \cdots  \\ 
\binom{n+r-1}{n+r-1} & \binom{n+r}{n+r-2} & \cdots  & \binom{n+2r-2}{n}%
\end{array}%
\right\vert _{r\times r}$\newline

$=\left\vert 
\begin{array}{ccccc}
1 & \frac{(n+1)n}{2!} & \frac{(n+2)(n+1)n(n-1)}{4!} & \cdots  & \frac{%
(n+r-1)(n+r-2)\cdots (n-r+2)}{2(r-1)!} \\ 
1 & \frac{(n+2)(n+1)}{2!} & \frac{(n+3)(n+2)(n+1)n}{4!} & \cdots  & \frac{%
(n-r)(n+r-1)\cdots (n-r+3)}{2(r-1)!} \\ 
1 & \frac{(n+3)(n+2)}{2!} & \frac{(n+4)(n+3)(n+2)(n+1)}{4!} & \cdots  & 
\frac{(n-r+1)(n+r-1)\cdots (n-r+4)}{2(r-1)!} \\ 
\cdots  & \cdots  & \cdots  & \cdots  & \cdots  \\ 
1 & \frac{(n+r)(n+r-1)}{2!} & \frac{(n+r+1)(n+r)(n+r-1)(n+r-2)}{4!} & \cdots 
& \frac{(n+2r-2)(n+2r-3)\cdots (n+1)}{2(r-1)!}%
\end{array}%
\right\vert _{r\times r}=\cdots .$

Set $6,7$ in above formula, we can get (\ref{cc3}) and (\ref{cc4}).

\end{document}